\documentclass[11pt]{amsart}
\usepackage{amsmath,amsthm,amscd}
\usepackage{amssymb,amsmath}
\usepackage{color}
\usepackage[all,cmtip]{xy}

\renewcommand{\mod}{\mathrm{mod\;}}
\newcommand{\Ext}{\mathrm{Ext}}

\newcommand{\rad}{\mathrm{rad}}
\newcommand{\hd}{\mathrm{hd}}
\newcommand{\A}{\mathbf{\rm A}}
\newcommand{\D}{\mathbf{\rm D}}
\newcommand{\E}{\mathbf{\rm E}}
\newcommand{\TL}{\mathbf{\rm L}}
\newcommand{\DL}{\mathbf{\rm DL}}

\newtheorem{Theorem}{Theorem}[section]
\newtheorem{Lemma}[Theorem]{Lemma}
\newtheorem{Proposition}[Theorem]{Proposition}
\newtheorem{Corollary}[Theorem]{Corollary}
\newtheorem{Definition}[Theorem]{Definition}
\newtheorem{Remark}[Theorem]{Remark}

\begin{document}

\parskip5pt
\parindent0pt

\title{Chebyshev polynomials on symmetric matrices}


\author{Karin Erdmann}
\address{Mathematical Institute, University of Oxford, UK}
 \email{erdmann@maths.ox.ac.uk}

\author{Sibylle Schroll}
\address{Department of Mathematics, University of Leicester, UK}
\email{schroll@mcs.le.ac.uk}

\thanks{ 
This work was supported by the Engineering and Physical Sciences Research Council [grant number EP/D077656/1] as well as by the Leverhulme
Trust through an Early Career Fellowship.}

\begin{abstract}
In this paper we evaluate  Chebyshev polynomials of the second-kind on a class of
symmetric integer matrices, namely on adjacency matrices of simply
laced Dynkin and extended Dynkin diagrams. As an application of these results we
explicitly calculate minimal projective resolutions of simple
modules of symmetric algebras with radical cube zero that are of
finite and tame representation type.
\end{abstract}

     \maketitle

\hsize440pt\hoffset-30pt

\baselineskip=16pt

 \section{Introduction}

Chebyshev polynomials are a sequence of recursively defined
polynomials. They appear in many  areas of mathematics such as
numerical analysis, differential equations, number theory and
algebra \cite{R}. Although they have been known and studied for a
long time, they  continue to play an important role in recent
advances in many subjects, for example in numerical integration,
polynomial approximation, or spectral methods ({\it e.g.}
\cite{MH}). It is interesting to note that they also play an
important part in the representation theory of algebras ({\it e.g.}
\cite{C, EF, GL, RX}). There are several closely related Chebyshev
polynomials. Amongst these, the polynomials usually referred to as
Chebyshev polynomials of the first kind, and Chebyshev polynomials
of the second kind are the ones that often naturally appear; for
example, they arise as solutions  of special cases of the
Sturm-Liouville differential equation or in dimension counting in
representation theory ({\it e.g.}   \cite{BM, E, MZ}).

For Chebyshev polynomials of the first and the second kind, the recursive definition is equivalent to a definition by a determinant
formula.  Symmetric integer matrices are a key to this definition.
 Focusing on the polynomials of the second kind,
we exhibit some surprising properties of Chebyshev polynomials in relation to these symmetric matrices.
In fact the symmetric matrices we consider are adjacency matrices of Dynkin diagrams and extended Dynkin diagrams.
   Dynkin diagrams play an import role in Lie theory, where they give a classification
of root systems. However, they also appear in areas that have no
obvious connection to Lie theory as for example in singularity
theory where they are linked to Kleinian singularities, or, for
example, in representation theory of algebras where they classify
symmetric algebras of radical cube zero of finite and tame
representation type~\cite{Benson}.

The motivation for the study of the Chebyshev polynomials evaluated on
matrices comes from the representation theory of the algebras
classified in ~\cite{Benson}. We begin with a detailed study of the
Chebyshev polynomials evaluated on adjacency
matrices of Dynkin diagrams, where we show that in the case of
Dynkin diagrams, the families of polynomials are periodic and
in the case of the extended Dynkin diagrams the families grow
linearly. We then show as an application how the general results we
obtain can be applied to the representation theory of the symmetric
algebras of radical cube zero.
We will see that the Chebyshev polynomials govern the minimal projective resolutions for these algebras and that they give rise to a method to calculate the constituents of minimal projective resolutions of simple modules.



We will now outline the content of this paper. In the next section
we recall the definition of Chebyshev polynomials of the second kind,
we define the polynomials we will be working with and we introduce the Dynkin diagrams together with a labeling of these diagrams which we will use throughout the paper.
 In
section~\ref{Evaluation} we evaluate Chebyshev polynomials on the adjacency matrices of Dynkin
diagrams
 and extended Dynkin
diagrams.
In section~\ref{Application} a link with the representation theory
of symmetric algebras of radical cube zero of finite and tame
representation type is described. In particular, we show how the
results of section~\ref{Evaluation} can be used to calculate minimal
projective resolutions of the simple modules. Finally in
section~\ref{Symmetric} we show a more general result on Chebyshev
polynomials evaluated on positive symmetric matrices.



\bigskip



\section{Definitions}

\subsection{Chebyshev Polynomials}

We briefly recall the definition of Chebyshev polynomials of the
second kind; good references are~\cite{R, S}.

The Chebyshev polynomial of the second kind $U_n(x)$ of degree $n$ is defined by
$$U_n(x) = {\rm sin}(n+1)\theta/{\rm sin} \theta \;\;\;\ \mbox{where} \; x= {\rm cos} \theta.$$
From this definition the following recurrence relations with initial conditions can be deduced
$$ U_n(x) = 2x U_{n-1}(x) - U_{n-2}(x) \;\;\; \mbox{with} \; U_0(x)=1 \; \mbox{and} \;  U_1(x) = 2x. $$
Furthermore an  easy calculation shows (see also~\cite{S} page 26) that $U_n(x) = \det (2xI_n - A_n)$ where $A_n$ is
an $n \times n$ matrix that has zeros everywhere except directly above and directly below the diagonal where all the entries are
equal to one.

We will work with the version of the Chebyshev polynomial defined by
$$f_n(x) = \det (xI - A_n)$$
 where $A_n$ is defined  as above, and $I$ is the identity matrix.  These polynomials are also sometimes called Dickson
 polynomials of the second kind.
An easy calculation shows that $f_n(x)$ is defined by the
recurrence relation $f_n(x) = x f_{n-1}(x) - f_{n-2}(x)   $ with
initial conditions $f_0(x)=1$ and $f_1(x) = x$. Furthermore, we set
$f_{-1}(x)=0$. All matrices have entries in  $\mathbb{Q}$.

\subsection{Dynkin diagrams and adjacency matrices}\label{Dynkin}

The Dynkin diagrams and extended Dynkin diagrams we are going to
consider are the ones of type $\A, \D, \E$ and $\widetilde{\A}$, $\widetilde{\D}$, $\widetilde{\E}$ as well as the diagrams of type
$\TL$, $\widetilde{\TL}$, and $\widetilde{\rm{DL}}$ (see Appendix).

Let $G$ be an undirected graph with $n$ vertices labeled by the set $\{1, 2, \ldots, n\}$.  The adjacency matrix of $G$ is a $n \times n$ matrix
where the entry in position $(i,j)$ is given by the number of edges between the vertices $i$ and $j$.
Modulo conjugation by a permutation matrix, the adjacency matrix is independent of the choice of labeling.

However, in what follows we work with particular  adjacency matrices
corresponding to a particular labeling of the graphs. We refer the
reader to the Appendix for the labeling of the diagrams we
have chosen. In the case of the Dynkin diagrams this labeling corresponds to the Dynkin labeling. It
follows from Lemma~\ref{trivfacts} in the next section that our
results are, up to permutation, independent of the labeling.

Our goal is to evaluate the Chebyshev polynomials of the second kind
on the adjacency matrices of the Dynkin diagrams and the extended
Dynkin diagrams as well as the diagrams of type $L$, $\widetilde{L}$, and $\widetilde{DL}$.

\section{Evaluating Chebyshev polynomials}\label{Evaluation}

\subsection{Evaluating Chebyshev polynomials on matrices - general results}
In this section we list some useful facts about evaluating the
Chebyshev polynomials on matrices. We begin with
a lemma collecting some straightforward facts.

\begin{Lemma}\label{trivfacts}
For any square matrix $M$, we have $Mf_k(M) = f_k(M)M$ and for any
symmetric matrix $S$, the matrix $f_k(S)$ is symmetric for all $k
\geq 0$.  If $T$ is an invertible $n \times n$ matrix then for all
$n \times n$ matrices $M$ we have $ f_k(TMT^{-1}) = T f_k(M)
T^{-1}$.
\end{Lemma}


\bigskip

\begin{Definition}  Let $X$ be a square matrix, we say that the
sequence of matrices $(f_k(X))_{k \geq 0}$ is periodic of period
$\leq p$, if $p>1$  and if $p$ satisfies $f_{p-1}(X)= 0$ and
$f_p(X)=I$.
\end{Definition}

{\bf Remark:} Suppose that $(f_k(X))_{k \geq 0}$ is periodic of
period $\leq p$. Then for any integer $k$, we can write $k=qp+r$
with $0\leq r < p$, and $f_k(X) = f_r(X)$.

\bigskip

\begin{Lemma}\label{Periodic}  Assume $X$ is a square matrix such that $f_d(X)=0$ for some
$d>1$. Then for $0\leq k\leq d+1$, we have
 \[f_{d+k}(X) + f_{d-k}(X)
= 0. \leqno{(*)}
\]
Moreover, $f_{2d+1}(X)=0, \ f_{2d+2}(X)=I$ and hence
$f_{2d+3}(X)=X$. Therefore the sequence $(f_k(X))_{k \geq 0}$ is
periodic, of period $\leq 2d+2$.
\end{Lemma}

\bigskip

{\it Proof: } \ The recursion for the Chebyshev polynomials can be
rewritten as
$$xf_m(x) = f_{m+1}(x) + f_{m-1}(x), \ \ \mbox{for} \ \ m\geq 0.
$$
Consider now $x^kf_d(x)$ for $0\leq k\leq d+1$  and substitute
$x=X$. Since $f_d(X)=0$, induction on $k$ will show that $X^kf_d(X)
= f_{d+k}(X) + f_{d-k}(X)=0$, for $0\leq k\leq d+1$.

\bigskip

The case $k=0$ is clear. Assume now that the statement is true for
all $j$ where $j\leq k$ and suppose $k < d+1$. Then
\begin{align*} 0 = X^{k+1}f_d(X) =
& X[X^kf_d(X)]\cr =& X[f_{d+k}(X) + f_{d-k}(X)] \cr =& f_{d+k+1}(X)
+ f_{d+k-1}(X) + f_{d-k+1}(X) + f_{d-k-1}(X)\cr =&   f_{d+k+1}(X) +
X^{k-1}f_d(X) + f_{d-k-1}(X) \cr =& f_{d+k+1}(X) + 0 + f_{d-k-1}(X)
\end{align*}
as required.

For the last part, let $k=d+1$ and recall that $f_{-1}(X)=0$, hence
$f_{2d+1}(X)=0$. Then let $k=d$, and recall that $f_0(X)=I$ which
implies $f_{2d}(X)= -I$. Since $f_{2d+1}(X)=0$ we can  apply (*)
with $2d+1$ instead of $d$. We then obtain $f_{2d+2}(X) + f_{2d}(X)
= 0$ and hence $f_{2d+2}(X) = I$ and finally $f_{2d+3}(X) = X\cdot I
- 0 = X$. $\Box$

\bigskip

If for some $d\geq 1$ the matrices $f_{d}(X)$ and $f_{d+1}(X)$ are
equal, then  periodicity  follows by Lemma~\ref{Periodic}.

\medskip

\begin{Lemma}\label{Equality}
Assume $X$ is a square matrix such that $f_d(X) = f_{d+1}(X)$ for
some integer $d\geq 1$. Then for $1 \leq k \leq d+1$, we have
$$f_{d+1+k}(X) = f_{d-k}(X).
$$
In particular, $f_{2d+2}(X)=0.$
\end{Lemma}

\bigskip

{\it Proof:} Assume $f_d(X) = f_{d+1}(X)$. Then we have
$$f_{d-1}(X) = Xf_d(X) - f_{d+1}(X) \ \ \mbox{and} \ \  f_{d+2}(X) = Xf_{d+1}(X) - f_d(X)
$$
and hence $f_{d-1}(X) = f_{d+2}(X)$. For the inductive step, assume
the statement is true for $1\leq m\leq k$. Then
\begin{align*}
f_{d-k}(X) = & Xf_{d-k+1}(X) - f_{d-k+2}(X) \cr = & Xf_{d+k}(X) -
f_{d+k-1}(X)\cr  = & f_{d+k+1}(X). \ \ \ \ \ \;\;\;\;\;\; \Box
\end{align*}

Some matrices allow a reduction, based on the following lemma,  which
gives a criterion to determine when a sequence of matrices has
linear growth.

\begin{Lemma}\label{GeneralReduction}  Assume $X$ is a square matrix such
that for some matrix $Z$ and for some $q\geq 2$  we have
 $f_q(X)=f_{q-2}(X) + Z$, and where  $ZX=XZ=2X$. Then

(a) \ for $1\leq t\leq q-1$, we have
$$f_{q+t}(X) = \left\{\begin{array}{ll} f_{q-2-t}(X) + 2f_t(X) & \mbox{ $t$ odd} \cr
                                         f_{q-2-t}(X) + 2f_t(X) + (-1)^{\frac{t}{2}+1}(2I-Z) &  \mbox{ $t$
                                         even},
\end{array}
\right.
$$


(b) \ if $Z=2I$ ,  for $t\geq -1$,  we have $f_{2q+t}(X) =
2f_{q+t}(X) - f_t(X)$,

(c) \ if $Z=2I$ and $m= rq  + u$ where $-1\leq u\leq q-2$ and $r
\geq 2$, we have

$$f_{rq+u}(X) = rf_{q+u}(X) - (r-1)f_u(X).$$
\end{Lemma}

\bigskip

\begin{Corollary}\label{GeneralReductionCorollary}
 Assume $X$ is a square matrix
such that for some matrix $Z$ we have
 $f_q(X)=f_{q-2}(X) + Z$, for $q\geq 2$ and where  $ZX=XZ=2X$ and where $Z \neq 2I$. Then
$$f_{2q}(X) - f_{2q-2}(X) = \left\{\begin{array}{ll} 2Z -2I & \mbox{ $q$ odd} \cr
                                         2Z-2I + (-1)^{\frac{q-2}{2}+1}(4I-2Z) &  \mbox{ $q$ even.}
\end{array}
\right.
$$
\end{Corollary}

\bigskip

{\it Proof:} By the Chebyshev recursion formula we have $f_{2q}(X) =
X f_{2q-1}(X) - f_{2q-2}(X)$  and $2q-1 = q + (q-1)$ and
$2q-2 = q + (q-2)$. Therefore we can apply
Lemma~\ref{GeneralReduction} a) and the result follows.  $\Box$

\bigskip

{\bf Remark:} (I) \ \ Later, the matrix $Z$ will often be equal to
$2I$. In this case (a) is equal to $f_{q+t}(X) = f_{q-2-t}(X) +
2f_t$. Furthermore, an easy calculation shows that (c) becomes
$f_{rq+u}(X) = rf_{q-2-u}(X) + (r+1)f_u(X)$.

(II) \ \ Observe that part (a) describes $f_k(X)$ for $q < k \leq
2q-1$. Since every natural number $m \geq 2q-1$ can be written in
the form given in part (c), this gives a description of $f_k(X)$ for
all $k\geq q$. It thus follows from the formula in part (c) that a
sequence $(f_k(X))_k$ that satisfies the hypotheses of
Lemma~\ref{GeneralReduction} with $Z=2I$ has linear growth.

\bigskip

{\it Proof of Lemma~\ref{GeneralReduction}: } (a) \ We use induction
on $t$. For $t=1$, we have
\begin{align*}
f_{q+1}(X) = & Xf_{q}(X) - f_{q-1}(X)\cr =& Xf_{q-2}(X) + XZ -
f_{q-1}(X)\cr  = &Xf_{q-2}(X) + 2X - f_{q-1}(X) \cr =& f_{q-3}(X) +
2f_1(X)
\end{align*}
if we recall that $f_1(X)=X$.

For $t=2$, we have
\begin{align*} f_{q+2}(X) = Xf_{q+1}(X) - f_{q}(X) &= X(f_{q-3}(X) + 2f_1(X)) - f_{q}(X)\cr
&= Xf_{q-3}(X) + 2Xf_1(X) - (f_{q-2}(X) + Z)\cr &= f_{q-4}(X) +
2f_2(X) + 2I -Z.
\end{align*}

For the inductive step suppose first that $t$ is odd. Then we have
\begin{align*} f_{q+t+1}(X) = & Xf_{q+t}(X) - f_{q+t-1}(X)\cr
=& X[f_{q-2-t}(X) + 2f_{t}(X)] - [f_{q-2-(t-1)}(X) + 2f_{t-1}(X)
+(-1)^{\frac{t-1}{2}+1}(2I-Z) ]
 \cr
=& f_{q-2-(t+1)}(X) + 2f_{t+1}(X) +(-1)^{\frac{t+1}{2}+1}(2I-Z).
\end{align*}

Now suppose that $t$ is even. Then we have

\begin{align*} f_{q+t+1}(X) = & Xf_{q+t}(X) - f_{q+t-1}(X)\cr
=& X[f_{q-2-t}(X) + 2f_{t}(X) +(-1)^{\frac{t}{2}+ 1}(2I-Z) ] -
[f_{q-2-(t-1)}(X) + 2f_{t-1}(X) ]
 \cr
=& f_{q-2-(t+1)}(X) + 2f_{t+1}(X) + (-1)^{\frac{t}{2}+1}(2X-ZX) \cr
=& f_{q-2-(t+1)}(X) + 2f_{t+1}(X)
\end{align*}
since $ZX =2X$.


\bigskip

(b) \ The case $t=-1$ follows from  part (a). Let $t=0$, then
\begin{align*} f_{2q}(X) = & Xf_{2q-1}(X) - f_{2q-2}(X)
\cr = & X[2f_{q-1}(X)] - [2f_{q-2}(X) +f_0(X)] \cr = & 2f_{q}(X) - I
\end{align*}
where the equality  $f_{2q-2}(X)= f_0(X) +2f_{q-2}(X)$ follows from
part (a). Let $t\geq 1$, and assume the equation holds for  $t-1$
and $t-2$, then
\begin{align*}
f_{2q+t}(X) =& Xf_{2q+(t-1)}(X) - f_{2q+(t-2)}(X) \cr
 =& X[2f_{q+t-1}(X) -
f_{t-1}(X)] - [2f_{q+t-2}(X) - f_{t-2}(X)]\cr =& 2f_{q+t}(X) -
f_{t}(X).
\end{align*}

\bigskip

(c) \ The case $r=2$ follows from part (b). Assume now $r\geq 3$ and
write $rq+u = 2q+t$ where  $t=(r-2)q+u$. Then by part (b) we have
$$
f_{rq+u}(X) = f_{2q+t}(X) = 2f_{q+t}(X) - f_t(X). \leqno{(*)}$$

If $r=3$ then $2q+u=q+t$ and using part (b) again we have
$f_{q+t}(X) = 2f_{q+u}(X) - f_u(X)$. Substituting this into equation
(*) gives
$$f_{3q+u}(X) = 2[2f_{q+u}(X) - f_u(X)] - f_{q+u}(X) = 3f_{q+u}(X) - 2f_u(X)
$$
as required.

Now let $r>3$.  We write $(r-1)q+u=q+t$ and $(r-2)q+u=t$. Then
equation (*) gives
\begin{align*} f_{rq+u}(X)= &2f_{q+t}(X) - f_t(X)\cr =
& 2[f_{(r-1)q+u}(X)] - f_{(r-2)q+u}(X) \cr =& 2[(r-1)f_{q+u}(X) -
(r-2)f_u(X)] - [(r-2)f_{q+u}(X) - (r-3)f_u(X)]\cr =& rf_{q+u}(X) -
(r-1)f_u(X)
\end{align*}
and this completes the proof. $\Box$

\bigskip

\begin{Definition}\label{Zero}
\normalfont Given some $n \times n$ matrix $X$, we define $X^0$ to be
the matrix obtained from $X$ by reversing the entries in each row.
In particular, in $I^0$ the $i,n-i+1$-th entries are equal to $1$
for all $1 \leq i \leq n$  and all other entries equal to zero.
\end{Definition}

Note that $(X^0)^0 = X$, and that $X^0 = X\cdot I^0$. Furthermore,
we write $I^0_n$ if we need to specify the size of the matrix.

\bigskip

With this we have the following variation of
Lemma~\ref{GeneralReduction}. The proof is a straightforward
modification of that of Lemma~\ref{GeneralReduction} and we leave
the details to the reader.

\begin{Lemma}
\label{GeneralReductionZero} Assume $X$ is a square matrix such
 that for some matrix $Z$ and for some $q\geq 2$
 we have $f_q(X)=f_{q-2}(X) + Z$,  and where $XZ=2X^0$. Then

(a) \ for $1\leq t\leq q-1$, we have
$$f_{q+t}(X) = f_{q-2-t}(X) + 2f_t(X)^0.$$
In particular, $f_{2q-2}(X) = 2f_{q-2}(X)^0 + I$ and $f_{2q-1}(X) =
2f_{q-1}(X)^0$,

(b) \ for $t\geq -1$ we have $f_{2q+t}(X) = 2f_{q+t}(X)^0 - f_t(X)$,

(c) \ if  $m= rq  + u$ where $-1\leq u\leq q-2$ and $r\geq 2$, we
have
$$f_{rq+u}(X) = \left\{\begin{array}{ll} rf_{q+u}(X)^0 - (r-1)f_u(X)
 & r \mbox{ even } \cr
rf_{q+u}(X) - (r-1)f_u(X)^0 & r \mbox{ odd. }
\end{array}
\right.
$$
\end{Lemma}

{\bf Remark:} Note if $Z$ is equal to  $2I^0$ then (c) can be
rewritten as $$f_{rq+u}(X) = \left\{\begin{array}{ll}
rf_{q-u-2}(X)^0 + (r+1)f_u(X)
 & r \mbox{ even } \cr
rf_{q-u-2}(X) + (r+1)f_u(X)^0 & r \mbox{ odd. }\end{array}\right.$$

\bigskip

\subsection{Substituting type $A$}

In this paragraph we evaluate the Chebyshev polynomials on the
adjacency matrix of a Dynkin diagram of type $A$. This will
be the basis for the calculation of almost all the other finite
types as well as  all the extended types.

Assume $A$ is the adjacency matrix of a  Dynkin diagram of
$A_n$, $n \geq 2$. Then by the Cayley-Hamilton Theorem  we know that
$f_n(A)=0$ and by Lemma~\ref{Periodic} the sequence $(f_k(A))_{k
\geq 0}$ is periodic of period  $\leq 2n+2$.

\bigskip

Let  $\Theta_k$ be the subset of
$\Theta=\{ (i,j): 1\leq i\leq n, 1\leq j\leq n\}$ given by
$$\Theta_k = \{ (i,j): k+2\leq i+j\leq 2n-k, \ \ i+j\equiv k (\mod 2), -k\leq j-i\leq k\}.
$$

We think of this as a subset of $\mathbb{N}\times \mathbb{N}$.
Note that $\Theta_k$ consists of the points  in $\Theta$ of parity
$i+j\equiv k  (\mod 2) $ which lie in the rectangle given by the
lines
\begin{equation}\label{grid}x+y=k+2, \ x+y=2n-k, \ y-x=k, \ y-x= -k
\end{equation}

This rectangle has corners $(1,k+1), (k+1,1)$ and $(n-k, n), \ (n,
n-k)$. In particular, $\Theta_n = \emptyset$ and $\Theta_{-1} = \emptyset$.

Write  $E_{ij}$ for the usual matrix unit.

\begin{Proposition}\label{linalgA}
For $-1 \leq k \leq n $,
\begin{equation}\label{fomula1typeA} f_k(A) = \sum_{(i,j)\in \Theta_{k } }
E_{ij}.
\end{equation}
In particular, $f_n(A)=0$.
\end{Proposition}

\bigskip

\begin{Corollary}
The family $(f_k(A))_k, k\geq -1$ is periodic of period $\leq 2n+2$.
\end{Corollary}

{\it Proof:} It follows from Lemma~\ref{Periodic} that for $0 \leq k
\leq n+1$, $f_{n+k}(A) + f_{n-k}(A) = 0.$ and that the sequence
$(f_k(A))_{k \geq 0}$ is periodic of period at most $2n+2$.  $\Box$

Furthermore, $f_k(A) = 0$ if and only if $k = mn + m -1 = m (n+1)
-1$ for $m\geq 0$ and the entries of $f_k(A)$ are known, for all $k
\geq 0$.


\vspace{1cm}

{\it Proof of \ref{linalgA}: } Let $-1 \leq k \leq n$ and let $W_k$
be the expression on the right hand side  of (\ref{fomula1typeA}). To prove (\ref{fomula1typeA}) it is enough to
show that
$$AW_k = W_{k+1}+W_{k-1}.
$$
We have $(\sum_{l=1}^{n-1} E_{l,l+1})E_{ij} = E_{i-1,j}$, except in
the case  $i=1$ where $(\sum_{l=1}^{n-1} E_{l,l+1})E_{1j} = 0$.
Note, however, that this only occurs for $(1,k+1)$, that is in the
`top corner' of $\Theta_k$. Similarly,
$(\sum_{l=1}^{n-1}E_{l+1,l})E_{ij} = E_{i+1,j}$ except in the case
$i=n$ where $(\sum_{l=1}^{n-1}E_{l+1,l})E_{nj} = 0$. This occurs
only for $(n, n-k)$, that is the `bottom corner' of $\Theta_k$.
Therefore $AE_{ij} = E_{i-1,j} + E_{i+1, j}$ for all $(i,j)\in
\Theta_k$, with the two exceptions as described above.

\medskip

If we visualize the $E_{ij}$ occurring in $W_k$ as grid points in
the rectangle defined by $\Theta_k$ then the terms of $AW_k$ are
obtained  by replacing each $(i,j)$ in this rectangle by the two
points below and above, namely $(i-1,j)$ and $(i+1,j)$ (with the
exceptions of the top corner where $(i+1,j)$ is missing, and the
bottom corner where $(i-1,j)$ is missing).

Following this process we obtain only the points $(r,s)$ with
$r+s\equiv k -1$(mod 2) and each point $(r,s)$ inside the rectangle
defined by (\ref{grid}) is obtained twice. Additionally we get all
points $(r,s)$ lying on the lines
$$x+y=k+1, \ x+y = 2n-k+1, \mbox{ and } y-x=k+1, \ y-x=-k-1
$$
exactly once.

Thus we obtain precisely the points corresponding to $W_{k+1}$ and
to $W_{k-1}$, while the points in the intersection appear twice.
 $\Box$

\bigskip

\subsection{Substituting type $L$}\label{TypeL}

In this paragraph we evaluate the Chebyshev polynomials on the
adjacency matrix of a  diagram of type $L$. Let $L$ be
the adjacency matrix of a diagram of type ${\rm L}_n $.
Then we can express $L$ in terms of $A$ such that $L= A+E_{11}$
where $A$ is the matrix of type $A_n$ as in the previous section.
We will now  express
the matrices $f_k(L)$ in terms of the matrices
$f_k(A)$. For $k=1, 2, \ldots, n-1$ we define
$$T_k:= \sum_{1\leq i, j,  i+j\leq k+1} E_{ij}.
$$
That is, each entry in the upper left corner, up to the line $i+j=k+1$,  is equal to  $1$ and all other entries are zero.

\begin{Proposition}\label{SubL}

(a) \ For $k=1, 2, \ldots, n-1$ we have $f_k(L) = f_k(A) + T_k$.\\
(b) \ $f_n(L) = f_{n-1}(L).$\\
(c) \ For $1\leq k\leq n+1$ we have $f_{n-1+k}(L) = f_{n-1-k}(L)$.
In particular $f_{2n-1}(L)=I$ and
$f_{2n}(L)=0$.\\
(d) \ We have $f_{4n+1}(L)=0$ and  $f_{4n+2}(L)=I$.
\end{Proposition}

\begin{Corollary}
The sequence $(f_k(L))_k$ is periodic of period $\leq 4n+2$.
\end{Corollary}

\bigskip

{\it Proof } \ (a)
 The proof is by induction. Clearly, $f_1(L) = f_1(A) + T_1$. Suppose that  $k < n-1$ and that
$f_k(L) =  f_k(A) + T_k$  holds. Then
$$ \begin{array}{ccl}
f_{k+1}(L)  &=& L f_{k}(L) -f_{k-1}(L) \\ \\
&=&  (A +E_{11}) ( f_k(A) + \sum_{i=1}^{k} \sum_{j=1}^{k+1-i} E_{ij}) - f_{k-1}(L)  \\ \\
&=& Af_k(A) +  \sum_{i=1}^{k} \sum_{j=1}^{k+1-i} (E_{i+1j} + E_{i-1j}) + E_{11} f_k(A) +
 E_{11} \sum_{i=1}^{k} \sum_{j=1}^{k+1-i} E_{ij} \\ \\
&& - f_{k-1}(L)  \\ \\
&=&  Af_k(A) +  \sum_{i=1}^{k-1} \sum_{j=1}^{k-i} E_{ij} +  \sum_{i=2}^{k+1} \sum_{j=1}^{k+2-i} E_{ij}
+ \sum_{j=1}^{k+1} E_{1j}  - f_{k-1}(L) \\ \\
&=& Af_k(A) - f_{k-1}(A) + \sum_{i=1}^{k-1} \sum_{j=1}^{k-i} E_{ij} +  \sum_{i=1}^{k+1} \sum_{j=1}^{k+2-i} E_{ij}
-  \sum_{i=1}^{k-1} \sum_{j=1}^{k-i} E_{ij}\\ \\
&=& f_{k+1}(A) + \sum_{i=1}^{k+1} \sum_{j=1}^{k+2-i} E_{ij}.\\ \\
\end{array}.$$

(b) \ The calculation in part (a) also holds when $k=n-1$ and therefore
$$ \begin{array}{ccl}
f_n(L) &=&  L f_{n-1}(L) -f_{n-2}(L)  \\
&=& f_n(A) +  \sum_{i=1}^{n} \sum_{j=1}^{k+1-i} E_{ij} \\
&=& 0 + f_{n-1}(L_n). \\
\end{array}$$

Part (c) and (d) follow from Lemma~\ref{Equality}. $\Box$

\bigskip

\begin{Remark}\label{TkBk}\normalfont
(I) \ We keep the record that the calculation in part (a) shows that
for $k < n-1$
$$E_{11}T_k + AT_k - T_{k-1} = T_{k+1}.$$

(II) \ Later we have to use the matrices  of type $L$ but where the
labeling is reversed. Then we have the description as in the above
proposition, where the only change is that $T_k$ is replaced by
$B_k:= I^0T_k$. This matrix is obtained from $B_k$ by reflecting in
the line $i+j=n+1$. For this, the calculation in part (a) shows that
for $k<n-1$
$$E_{nn}B_k + AB_k - B_{k-1} = B_{k+1}.$$
\end{Remark}

\bigskip

\subsection{Substituting extended type L}

Fix $n\geq 2$ and let $\widetilde{L} = \widetilde{L}_{n}$ be the
adjacency matrix of a diagram of
 type ${\rm \widetilde{L}_{n}}$ as described in section~\ref{Dynkin}.
In this paragraph we evaluated the Chebyshev polynomial on
$\widetilde{L}$ and show that the family $f_k(\widetilde{L})$ is of
linear growth.

We note that $\widetilde{L}$ can be expressed in terms of the
adjacency matrix $A$ of $A_n$. Namely, $\widetilde{L} = A + E_{11} +
E_{nn}$. The next Proposition shows that the terms of
$f_k(\widetilde{L}_{n})$ are a sum of $f_k(A)$ and an upper and a
lower triangular matrix whose entries are all equal to 1. Recall the
definition of $T_k$ and $B_k$ from Section~\ref{TypeL}.

\begin{Proposition}\label{firsttermsEL}
For $1 \leq k \leq n-1$,  we have
$$ f_k(\widetilde{L}) = f_k (A) +
T_k + B_k.$$
\end{Proposition}



\bigskip

{\it Proof:} We proof the result by induction. The result holds  for
$k=1$ and a direct calculation shows that it also holds for $k=2$.
Suppose it holds for all $l \leq k$. Then by definition we have
$f_{k+1}(\widetilde{L}) = X f_k(\widetilde{L}) -
f_{k-1}(\widetilde{L})$ and by induction hypothesis this is equal to
\begin{eqnarray*} (A + E_{11} + E_{nn})(
f_k (A) + T_k + B_k) - f_{k-1} (A) - T_{k-1}-B_{k-1}
.\end{eqnarray*} To prove the stated formula we need
\begin{align*}
(E_{11}+E_{nn})(f_k(A) + T_k+B_k) + A(T_k+B_k) \cr - T_{k-1}-B_{k-1}
&= T_{k+1} + B_{k+1}
\end{align*}
This is true by the record kept in remark~\ref{TkBk}. $\Box$

Let $U$ be the $n \times n $
matrix all of whose entries are equal to 1. Recall also the
definition of $I^0$ and $X^0$, and note that
$2XI^0 = 2X^0$. As special cases of the previous, and
by applying Lemma~\ref{GeneralReductionZero},  we have therefore

\begin{Corollary}\label{extendedL}
We have that $f_{n-2}(\widetilde{L}) = U-I^0$, and
$f_{n-1}(\widetilde{L}) = U$, and $f_n(\widetilde{L}) = U+I^0$.
Hence
$$f_n(\widetilde{L})  = f_{n-2}(\widetilde{L}) + 2I^0$$
and the sequence $(f_k(\widetilde{L}))_{k\geq 0}$ has linear growth.
\end{Corollary}

\bigskip

\subsection{General setup for the remaining infinite families}\label{GeneralSetup}

All remaining infinite families of Dynkin and extended  Dynkin types are
based on the above calculations and fit into a more general set-up
described in this section. Namely, in this section we substitute a
symmetric block matrix $X$ of the form
$$X = \left(\begin{matrix}0&S\cr S^t & W\end{matrix}\right),
$$
where we assume that $SS^tS=2S$, and hence that  $S^tS$ and $SS^t$
are projections. Assume further that $S$ has rank one. Then for any
matrix $M$ of the appropriate size, $SMS^t$ is a scalar multiple of
$SS^t$.


The following lemma is a straightforward calculation.

\begin{Lemma}\label{recursionHSW}
Write $f_k(X)$ as  block matrix,
$f_k(X) = \left(\begin{smallmatrix} H_k & S_k \cr (S_k)^t & W_k\end{smallmatrix}\right).
$ Then we have
\begin{align*}H_{k+1} = & SS_k^t - H_{k-1}\cr
               S_{k+1} = & SW_k - S_{k-1} \cr
               W_{k+1} = & S^tS_k + WW_k - W_{k-1}.
\end{align*}
\end{Lemma}


We can express the adjacency matrices of the remaining infinite
families in question in terms of $S$ and
$W$ and we will apply Lemma~\ref{recursionHSW} as follows. Let
$\varepsilon_1$ be the row vector of length $n$ whose first entry is
equal to 1 and whose other entries are all equal to 0.

\noindent (1) \ For type $D_m$, $m\geq 3$, of size $m$,  the
adjacency matrix of type $D_m$ is given by choosing $W$ to be of
type $A_n$ with $n = m-2$, and $S$ to be the matrix with two rows,
each row equal to $\varepsilon_1$.


\noindent (2) \ For type $\widetilde{D}_{m+1}$, $m \geq 4$, of size
$m+2$, we take $S$ as in (1), and for $W$ we take the matrix $D$ of
type $D_{m}$ but with the labeling reversed, we will call this
matrix $V$, that is $V = I^{0} D I^{0}$.


\noindent (3) \ For type $\widetilde{A}_n$, $n \geq 3$, of size
$n+1$, we take $W=A_n$ and $S$ to be the matrix with one row equal
to $\varepsilon_1$.

\noindent (4) \ For type $\widetilde{DL}_m$, $m \geq 3$, of size
$m$, we take $W = I^{0} L_n I^{0}$ with $n=m-2$, that is $W$ is
$L_n$ but with reversed order, and we take $S$ as in (1) so that $S$
is a two row matrix with both rows equal to $\varepsilon_1$.

\bigskip

Take $X$ as above, and $S_k, W_k$ and $H_k$ as in the recursion in
Lemma~\ref{recursionHSW}.
We will now calculate the first few terms explicitly.

\medskip

\begin{Definition}\label{Invariants} We define invariants $c$ and $\lambda_c$  of $X$, to be the
first integer $c>1$, and the scalar $\lambda_c$,
 such that $Sf_c(W)S^t = \lambda_cSS^t$ is non-zero.
\end{Definition}

\bigskip

\begin{Proposition}\label{FormulasHSW} Let $c$ be as above and let $1\leq k\leq c+2$.  Then\\
(a) \
$$S_k = S(\sum_{i\geq 0} f_{k-1-2i}(W)), \ \  W_k = f_k(W) + \sum_{i\geq 0} \psi_{k-2-2i}
$$
where $\psi_x = \sum_{0\leq r\leq x} f_r(W)S^tS f_{x-r}(W)$.\\
(b) \ For $k < c+2$
$$H_k = \left\{\begin{array}{ll} 0 & k \ \mbox{odd}\cr
                                 I & k\equiv 0 \mbox{ mod 4}\cr
                                 SS^t-I & k\equiv 2 \mbox{ mod 4}.
\end{array}
\right.
$$
(c) \
$$H_{c+2} = \left\{\begin{array}{ll} (1+\lambda_c) SS^t - H_c & c \mbox{ even } \cr
 \lambda_c SS^t &   c \mbox{ odd}.
\end{array}
\right. $$
\end{Proposition}

We make the convention that $\psi_x=0$ and $f_x(W)=0$ if $x<0$.

\bigskip

{\it Proof: } (a) \ Induction on $k$. The cases $k=1$ and $k=2$ are
clear. So assume true for all $j$ with $1\leq j \leq k$, and suppose
$k < c+2$, then
$$ S_{k+1} = SW_k - S_{k-1}
=  S(f_k(W) + \sum_{i\geq 0} \psi_{k-2-2i}) - S(\sum_{j\geq 0}
f_{k-2-2j}(W)).
$$
Let $x \leq k-2$, since $k< c+2$ we have $x<c$, and therefore
$Sf_{x-y}(W)S^t=0$ for all $y\geq 0$. This implies that almost all
terms of $S\psi_x$ are zero, and
$$S\psi_x = SS^tSf_x(W) = 2Sf_x(W).
$$
Substituting this gives
$$S_{k+1} = Sf_k(W) + 2S(\sum_{i\geq 0} f_{k-2-2i}(W)) - S(\sum_{j\geq 0} f_{k-2-2j}(W)
$$
which proves the claim.

\medskip

Next, consider $W_{k+1}$, substituting the terms using the induction
hypothesis we get
\begin{align*} W_{k+1} =
 S^tS[\sum_{i\geq 0} f_{k-1-2i}(W)] & + Wf_k(W)  + \sum_{i\geq 0} W\psi_{k-2-2i} \cr
& - f_{k-1}(W) - \sum_{i\geq 0} \psi_{k-3-2i}.
\end{align*}

By the recursion, $Wf_k(W) - f_{k-1}(W) = f_{k+1}(W)$. Moreover, it
also follows from the recursion that, for $x\geq 0$,
$$W\psi_x - \psi_{x-1} = \psi_{x+1} - S^tSf_{x+1}(W).
$$
Substituting these, and noting that the terms $-S^tSf_{x+1}(W)$
cancel for all $x=k-2-2i$ gives the claim.

\medskip

(b)  and (c) \ The cases  $k=1$ and $k=2$ are clear. For the
inductive step, we have if $k< c+2$ that
$$H_{k+1} = S(S_k^t) - H_{k-1} =  S(\sum_{i\geq 0} f_{k-1-2i}(W))S^t -H_{k-1}.$$
With the assumption, $k-1-2i\leq c$ for all $i\geq 0$, and the only
way to get this equal to $c$ is for $k=c+1$ and $i=0$. Hence
$$S\sum_{i\geq 0} f_{k-1-2i}(W)S^t = \left\{\begin{array}{ll} Sf_c(W)S^t &  k-1=c \mbox{ odd}
\cr Sf_c(W)S^t + 2SS^t & k-1=c \mbox{ even} \cr SS^t & k-1 < c, k-1
\mbox{ even}\cr 0 & \mbox{else}
\end{array}
\right.
$$
and recalling that $Sf_c(W)S^t = \lambda_cSS^t$, this gives the
claim. $\Box$

\bigskip

\subsection{Substituting  type $D$}\label{Type D}

In this paragraph we evaluate the Chebyshev polynomials on the
adjacency matrices of  Dynkin diagrams of type $D$. Let $X$ be the
matrix associated to type $D_m$, such that $S$ and $W= A_n$ for
$m=n+2$ are as described in section~\ref{GeneralSetup}. According to
Proposition~\ref{linalgA} and the remark on $\Theta_k$ preceding it,
 the parameters defined in Definition~\ref{Invariants} are $c=2n$ and $\lambda_c=1$.
 Therefore we obtain the expressions of
$f_k(D_m)$ for $k\leq 2n+2$ from
Proposition~\ref{FormulasHSW}.

\bigskip

\begin{Lemma}\label{TypeD1} For $1\leq k< 2n+1$,  $f_k(D_m)\neq 0$   and $f_{2n+1}(D_m)=0$.
\end{Lemma}

\bigskip

{\it Proof: } \  First, we observe that

\begin{equation}\label{Sumf_kW} \sum_{i\geq 0} f_{2n-2i}(W)=0, \ \ \sum_{i\geq 0}
f_{2n-1-2i}(W) = 0.
\end{equation}

Namely, each of these  is a sum of terms of the form $f_{n+t}(W) +
f_{n-t}(W)$ for some $0\leq t\leq n+1$,
and by Lemma~\ref{Periodic} this sum is equal to zero.
Furthermore, if $1\leq r < 2n-1$ then the first row of $\sum_{i\geq
0} f_{r-2i}(W)$ is non-zero: for example, the entries of the first
row of $f_0(W)$ (or $f_1(W)$) do not cancel.

\bigskip

 We start by showing that for $1\leq k< 2n+1$,  $f_k(D_m)\neq 0$. Suppose $1\leq k < 2n+1$. Then by Proposition~\ref{FormulasHSW} we have
$$S_k = S[\sum_{i\geq 0} f_{k-1-2i}(W)].
$$
This is the two-row matrix where both rows are equal to the first
row of $\sum_{i\geq 0} f_{k-1-2i}(W)$. By the above this is non-zero
except when $k=2n$.

Thus for $k<2n$ we have $S_k\neq 0$ and therefore $f_k(D_m)\neq 0$.
If $k=2n$, then $k$ is even and by Proposition~\ref{FormulasHSW}$
H_k\neq 0$, and so  $f_{2n}(D_{2m})\neq 0$.

 \ Finally consider $f_{2n+1}(D_m)$. It follows from (\ref{Sumf_kW}) that $S_{2n+1}=0$,
and  Proposition~\ref{FormulasHSW} implies $H_{2n+1}=0$. Recall that
$f_{2n+1}(W)=0$, and thus
\begin{equation}W_{2n+1} = \sum_{i\geq 0} \psi_{2n-1-2i} \ = \ \sum f_r(W)SS^t f_s(W)
\end{equation}
 where the sum is over all $r, s\geq 0$ with $r+s\leq
2n-1$ and $r+s$ odd.

Given such $r, s$, define $r'$  and $s'$ by $r+r'=2n$ and $s'+s=2n$.
Then  for $0\leq k\leq n+1$, $r= n-k$ implies $r'=n+k$ and similarly
$s= n-k$ implies $s'=n+k$ or vice versa. This implies that $f_r(W) +
f_{r'}(W) =0 $ and  $f_s(W) + f_{s'}(W)=0$.

\medskip

It is clear that both $r+s'$ and $s+r'$ are odd. We now claim that
precisely one of $r+s'$ and $r'+s$ is strictly less than $2n$.

Assume for a contradiction that $r'+s \geq 2n$ and $r+s' \geq 2n$.
Then because both expressions are odd, they are both strictly larger
than $2n$ and we have $r'+s > 2n= r'+r$ and thus $s>r$. But $r+s' >
2n = s+s'$ implies $r>s$, a contradiction.

We get a similar contradiction if we assume $r'+s<2n$ and $r+s'< 2n$
and therefore exactly one of $r+s'$ and $r'+s$ is strictly less than
$2n$.

\bigskip

Suppose now that $r'+s< 2n$. Then the terms in (\ref{Sumf_kW}) where
labels of the form $r, r', s, s'$  occur are precisely
$$f_r(W)S^tSf_s(W) + f_{r'}(W)S^tSf_s(W).
$$
But this expression is zero since $f_r(W) + f_{r'}(W)=0$. $\Box$

\bigskip

The following two lemmas give more precise information about
particular entries of $f_k(D_m)$.

\bigskip

\begin{Lemma}\label{f_2nD_m}We have
$$f_{2n}(D_m) = \left\{ \begin{array}{ll} I & \mbox{ n even}
\cr & \cr
\left(\begin{matrix} I_2^0 & 0 \cr 0 & I_n
\end{matrix}\right)
 & \mbox{ n odd. }
\end{array}
\right.
$$
\end{Lemma}

\bigskip

{\it Proof: } Consider $W_{2n}= f_{2n}(W) + \sum_i \psi_{2n-2-2i}$.
We know $f_{2n}(W)=I$ and so we must show that $\sum \psi_{2n-2-2i}
= 0$. However this follows directly from an argument similar to the
one in the
 previous Lemma. Similarly one
shows that $S_{2n}=0$. The result then follows from Proposition~\ref{FormulasHSW}. $\Box$

\bigskip

\begin{Lemma}  Assume $1\leq k\leq 2n$. Then the last row of
$f_k(D_m)$ is equal to
$$\begin{array}{ll} \left(\begin{matrix} 0 & 0 & \varepsilon_{n-k}& &\end{matrix}\right) & 1\leq k < n \cr
                     \left(\begin{matrix} 1 & 1 & 0&...& 0 \end{matrix}\right) & k=n \cr
                    \left(\begin{matrix} 0 & 0 &  \varepsilon_{k-n}& &\end{matrix}\right)& n<k\leq
                    2n.
\end{array}
$$
\end{Lemma}

\bigskip

Hence the last row of $f_k(D_m)$ reversed is
$$\begin{array}{ll} \left(\begin{matrix}\varepsilon_{k+1}& & & &0&0\end{matrix} \right)& 1\leq k < n \cr
                    \left(\begin{matrix} 0&...&0 &  &1 &  1\end{matrix}\right)  & k=n \cr
                    \left(\begin{matrix} \varepsilon_{n-(k-n)+1} & 0&0 \end{matrix}\right)& n<k\leq
                    2n.
\end{array}
$$

\bigskip

{\it Proof: } (1) We need the last column of $S_k$, the transpose of
this gives the first two entries for the required last row.

Recall $S_k = S[\sum_{i} f_{k-1-2i}(W)]$. This is the 2-row matrix
where each row is equal to the first row of $\sum_i f_{k-1-2i}(W)$.
We need only the $1n$ entry of this sum.

The $1n$ entry of $f_x(W)$ is $1$ for $x=n-1$ and $-1$ for $x=n+1$
and is zero for any other $x\leq 2n$.

\bigskip

{\sc Case 1 } Assume $n$ is odd. Then $f_x(W)_{1n}=0$ for all odd
$x\leq 2n$. Set $x=k-1-2i$, so for $k$ even, the last column of
$S_k$ is zero.

Now let $k$ be odd, and consider the $x=k-1-2i$. We have
$f_x(W)_{1n}=1$ for $x=n-1$ and $=-1$ for $x=n+1$ and is zero
otherwise. It follows that the 1n entry of $\sum_{i} f_{k-1-2i}(W)$
in this case is equal to $1$ if $k-1=n-1$ and is zero otherwise, by
cancelation.

This shows that  the last column of $S_k$ is zero unless $k=n$ and
then it is of the form  ${1\choose 1}$.

\medskip

{\sc Case 2 } Assume $n$ is even. Then for $k$ odd (all $x$ even),
as before the last column of $S_k$ is zero. Assume $k$ is even,
consider $x=k-1-2i$. We have $f_x(W)_{1n}$ as before. It follows
that the $1n$ entry of $\sum_i f_{k-1-2i}(W)$ is equal to $1$ if
$k-1=n-1$ and zero otherwise. Again, the last column of $S_k$ is
zero unless $k=n$ and then it is of the form ${1\choose 1}$.

\bigskip

(2) Now  we determine the last row of $W_k$, recall from
Proposition~\ref{FormulasHSW} that $W_k = f_k(W) + \sum_{i\geq 0}
\psi_{k-2-2i}$.

\ \ Assume first that $1\leq k\leq n$. We claim that then  the last
row of $\sum_i \psi_{k-2-2i}$ is zero, hence the last row of $W_k$
is equal to $\varepsilon_{n-k}$ for $k<n$,  and is zero for
$k=n$.

Consider $(Sf_a(W))^t(Sf_b(W))$. If this has last row non-zero then
we must have that the first row of $f_a(W)$ has non-zero $1n$ entry.
This occurs only for $a=n-1$ or $a=n+1$ but here we have only
$a+b\leq k-2\leq n-2$. So this has last row equal to zero. This
implies the claim.

\ \ Now consider $k=n+r$ where $1\leq r\leq n$.  Then $f_{n+r}(W) =
-f_{n-r}(W)$ and this has last row equal to $-\varepsilon_r$. We
claim that $\sum_i \psi_{n+r-2-2i} = 2\varepsilon_r$. [This will
imply the statement.]

\medskip

We use  induction on $r$. Assume first that  $r=1$.

Then $\sum_i \psi_{n-1-2i} = \psi_{n -1} + \sum_{i>0}
\psi_{n-1-2i}$. In the sum, the last row is zero (by the argument in
the previous case). The last row of $\psi_{n-1}$ has non-zero
contribution only from $f_{n-1}(W)S^t$ and this is $2\varepsilon_1$

\medskip

For the inductive step, write
$$\sum_i \psi_{n+r-2-2i} = \psi_{n+r-2} + \sum _{i\geq 0}
\psi_{n+r-4-2i}.
$$
By the inductive hypothesis the sum is equal to
$2\varepsilon_{r-2}$. Now consider $\psi_{n+r-2}$. This has only two
terms with non-zero last row, and the last row of $\psi_{n+r-2}$ is
equal to the last row of
$$(Sf_{n-1}(W))^tSf_{r-1} + (Sf_{n+1}(W))^tSf_{r-3}.
$$
This is equal to $2\varepsilon_r - 2\varepsilon_{r-2}$. In total we
get the stated answer. $\Box$

\bigskip

Using 4.3 with $d=2n+1$,
the previous three Lemmas imply that
$f_{2m-3-k}(D_m) = f_{2m-3+k}(D_m)$ for $0\leq k\leq 2m-3$;   and the periodicity of $
(f_k(D_m))_{k\geq 0}$ follows.

\begin{Corollary}
The family $(f_k(D_m))$ is periodic of period $\leq 2(2m-2)$.
\end{Corollary}

\bigskip

\subsection{Substituting type $\widetilde{D}$}

Let $X$ be of type $\widetilde{D}_{m+1}$ of size $m+2 = n+4$. That is,
$X=\left(\begin{matrix}0 & S\cr S^t & V\end{matrix}\right)$, where
$V$ is equal to $D_{n+2}$ with
reversed order (explicitly $V = I^0D_{n+2}I^0$), and where
 $S = \left(\begin{matrix}1 & 0& \ldots & 0\cr
                         1 & 0 & \ldots & 0\end{matrix}\right)$ of size $2\times m$.  .
Recall from ~\ref{Invariants} the definition of the
invariants $c$ and $\lambda_c$.

\begin{Lemma} For $X$ as described above, we have $c=2n$ and $\lambda_c=1$.
\end{Lemma}

\bigskip

{\it Proof: } The parameter $\lambda_c$ is the $11$ entry of
$f_c(V)$, which is equal to the $nn$ entry of $f_c(D_m)$ when
this is non-zero the first time. The $nn$ entry of $f_k(D_m)$ is  the
$nn$ entry of the matrix $W_k$ occurring
in the recursion for type $D_m$,  with
the notation of section~\ref{Type D}.

Let $k\leq 2n$. We know from Proposition~\ref{FormulasHSW} that $W_k
= f_k(W) + \sum_{i\geq 0} \psi_{k-2-2i}$. We notice that
$(f_k(W))_{nn}\neq 0$ only if $f_k(W)= \pm I$ and this occurs for
the first time when $k=2n$.

\medskip

Now consider the $nn$ entry of $\psi_x$ for $x\leq 2n-2$. For this
we need the $nn$ entry of $(Sf_r(W))^t(Sf_s(W)$ for $r+s=x\leq
2n-2$. This is equal to
$$2(f_r(W))_{1n}\cdot f_s(W)_{1n}.
\leqno{(*)}
$$
The $1n$ entry of any $f_t(W)$ for $t \leq 2n-2$ is only non-zero
for $t=n-1$ or $n+1$. So the number (*) is zero for $r+s < 2n-2$. If
$r+s=2n-2$ then it is only non-zero for $r=s=n-1$ and then it is
equal to 1. $\Box$

\bigskip

By Proposition~\ref{FormulasHSW} we now obtain the expressions for
the matrices constituting $f_k(X)$ for $1\leq k\leq 2n+2$, and this
is enough to prove linear growth. Recall from Definition~\ref{Zero}
the definition of $I_2^0$.

\begin{Proposition}\label{TowardsTildeDm+1} We have
$$f_{2n+2}(X) - f_{2n}(X) = \left\{ \begin{array}{ll} 2I & n \mbox{ odd } \cr
 2\widetilde{I} & n \mbox{ even }
\end{array}
\right.
$$
where
$$\widetilde{I} =
\left(\begin{matrix}I_2^0 & 0 &0 \cr 0 & I_n & 0 \cr 0&0& I_2^0
\end{matrix}\right)
$$
Hence the sequence $(f_k(X))_k$ has linear growth.
\end{Proposition}

{\it Proof:} By \ref{FormulasHSW}, using $c=2n$ and $\lambda_c=1$,
 we have $H_{2n+2} - H_{2n} = 2SS^t-2H_{2n}
= 2I_2$ for $n$ odd, and $H_{2n+2} - H_{2n} = 2(SS^t-I_2)$ if $n$ is
even.

Now consider $S_{2n+2} - S_{2n}$.  Here most of terms of the sum
cancel and we obtain
$$ S_{2n+2} - S_{2n} =Sf_{2n+1}(V) = S\cdot 0 = 0.
$$

Finally, for  $V_{2n+2} - V_{2n}$ most of the terms of the sum
cancel as well and we are left with  $V_{2n+2} - V_{2n} =
f_{2n+2}(V) - f_{2n}(V) + \psi_{2n}$. Furthermore, we have
$$ f_{2n+2}(V) -
f_{2n}(V) = -2f_{2n}(V) =  \left \{
\begin{array}{ll}
-2I_{n+2} & \mbox{ $n$ even} \\
 -2U &  \mbox{ $n$ odd} \\
 \end{array} \right.$$
 where $U$ is the matrix $f_{2n}(D_m)$ reversed.

 It remains to calculate
$\psi_{2n}$.
 For that we note that
\begin{align*} (Sf_0(V))^t(Sf_{2n}(V)) = & 2E_{11} = (Sf_{2n}(V))^t(Sf_0(V))\cr
(Sf_1(V))^t(Sf_{2n-1}(V)) = & 2E_{22} = (Sf_{2n-1}(V))^t(Sf_1(V))\cr
\ldots & \ldots \cr
(Sf_{n-1}(V))^t(Sf_{n+1}(V)) = & 2E_{nn} =
(Sf_{n+1}(V))^t(Sf_{n-1}(V))\cr
(Sf_n(V))^t(Sf_n(V)) = &
2 (E_{n+1,n+1} + E_{n+1,n+2} + E_{n+2,n+1} + E_{n+2,n+2}).
\end{align*}
Thus $\psi_{2n} = \left(\begin{matrix} 4I_n & 0 \cr 0 & 2SS^t
\end{matrix}\right)$, and
$$V_{2n+2} - V_{2n} = \left\{\begin{array}{ll} 2I_{n+2} & n \mbox{ odd }\cr
 & \cr
\left(\begin{matrix} 2I_n & 0 \cr
                 0 & 2(SS^t-I_2)\end{matrix}\right)
& n \mbox{ even }
\end{array}
\right.
$$
The conclusion follows from Lemma~\ref{GeneralReduction} with $Z=2I$
for $n$ odd.
Assume $n$ is even, then applying Corollary~\ref{GeneralReductionCorollary} twice gives
$f_{8n+8}(X) - f_{8n+6}(X) = 2I$, and linear growth follows from Lemma~\ref{GeneralReduction}.
 $\Box$

\bigskip

\subsection{Substituting extended type $\widetilde{A}$}

Let now $X$ be the matrix of type $\widetilde{A}_{n}$, that is $X$
is an $n+1\times n+1$ matrix of the form $X=\left(\begin{matrix}0 &
S\cr S^t & W\end{matrix}\right)$, with
$$W=A_n \ \ \mbox{and} \ \ S = \left(\begin{matrix}1 & 0&  \ldots & 0&
1\end{matrix}\right).
$$
where  $S$ is a $1\times n$ matrix. One checks that the
invariants $c$ and $\lambda_c$ of $X$ as defined in
\ref{Invariants} are as follows.

\medskip

\begin{Lemma} Let $1\leq c$ be  minimal such that
$Sf_c(W)S^t = \lambda_cSS^t \neq 0$. Then $c=n-1$, and  $\lambda_c=1$.
\end{Lemma}

\bigskip

Therefore,  by Proposition~\ref{FormulasHSW} we get the matrices
$f_k(X)$ for $1\leq k\leq n+1$.

\begin{Proposition} We have
$$f_{n+1}(X) - f_{n-1}(X) = 2I.$$
Hence the sequence $(f_k(X))_k$ has linear growth.
\end{Proposition}

{\it Proof:} \ We have $H_{n+1}- H_{n-1} = 2$. Furthermore, $S_{n+1}
- S_{n-1} = Sf_n(W) =0.$ Consider $W_{n+1} - W_{n-1}$, this is equal
to
$$f_{n+1}(W) - f_{n-1}(W) + \psi_{n-1}$$
(since most of the terms of the sums cancel).

We know from type $A$ that $f_{n+1}(W) - f_{n-1}(W) = -2I^0$. A straightforward
calculation shows that
$$\psi_{n-1} = 2(I+I^0).$$
Combining these gives the  first  statement.
We get linear growth from Lemma \ref{GeneralReduction}.
$\Box$

\subsection{Substituting  type $\widetilde{\DL}$}

Let $X$ be the adjacency matrix of type $\widetilde{DL}_{m}$ with $m$
vertices, and $m=n+2$. With the notation of
section~\ref{GeneralSetup} $X$ has the following blocks:  $S =
\left(\begin{smallmatrix}1&0&0&\ldots &0\cr
                         1&0&0&\ldots &0\end{smallmatrix}\right)$, and $W=V$
where $V = I^0 L I^0$ is the adjacency matrix of  type $L_n$
reversed.

Recall that $f_{2n}(V)=0$ and $f_k(V)\neq 0$ for $1\leq k < 2n$.
Furthermore $f_{2n-1}(V) = I$. The invariants $c$  and $\lambda_c$
of $X$, given by $\lambda_cSS^t= Sf_c(V)S^t \neq 0$, and
$Sf_k(V)S^t=0$ for $1< k < c$ are as follows

\begin{Lemma} We have $c=2n-1$ and $\lambda_c=1$.
\end{Lemma}

{\it Proof} The parameter $Sf_k(V)S^t$ is the $(1,1)$ entry of
$f_k(V)$ which is equal to the $(n,n)$ entry of $f_k(L_n)$. It
follows from Section~\ref{TypeL} that this entry is equal to $1$ for
$c=2n-1$ and equal to $0$ for $1< k< 2n-1$. $\Box$

We can now apply Proposition~\ref{FormulasHSW}:

\begin{Proposition}  We have
$$f_{2n+1}(X) - f_{2n-1}(X) =  Z$$
where
$$Z = \left(\begin{matrix} SS^t & 0 \cr
0 & 2I_n\end{matrix} \right)
$$
and $XZ =ZX= 2X$. Hence the sequence $(f_k(X))_k$ has linear growth.
\end{Proposition}

{\it Proof: }  \ Using the formulae from
Proposition~\ref{FormulasHSW} with $c=2n-1$, we have $H_{2n+1} -
H_{2n-1} = Sf_{2n-1}(V)S^t = SS^t$. Next, consider $S_{2n+1} -
S_{2n-1}$. Here most of the sum cancels and the only expression
remaining is
$$Sf_{2n}(V) = S\cdot 0 = 0.
$$
Consider now $V_{2n+1} - V_{2n-1}$. This is equal to
$$f_{2n+1}(V) - f_{2n-1}(V) + \psi_{2n-1}.$$
By Proposition~\ref{SubL} {\it c)} combined with
Lemma~\ref{Periodic}, we have $f_{2n+1}(V) - f_{2n-1}(V) = -2I$. We
now calculate $\psi_{2n-1}$, this is equal to
$$\psi_{2n-1}=2[\sum_{r=0}^{n-1} (Sf_r)^t\cdot (Sf_{2n-1-r})].
$$
Similarly to the proof of Proposition 9.2 we find that
$(Sf_r)^t\cdot (Sf_{2n-1-r}) = 2E_{r+1,r+1}$, for $0\leq r\leq n-1$.
Hence $V_{2n+1} - V_{2n-1} = -2I+4I = 2I$, as required. This proves
that the difference $f_{2n+1}(X) - f_{2n-1}(X)$ is as stated, and
one checks that $XZ = 2X$. Now by applying
Corollary~\ref{GeneralReductionCorollary} first with $q' = 2q$ we
obtain $f_{q'}(X) - f_{q'-2}(X) =  Z'$ where $Z' =
\left(\begin{matrix} I_2^0 & 0 \cr 0 & 2I_n\end{matrix} \right)$ and
where $Z' X= XZ' = 2X$. Then we set $q^{''} =  2q' =4q$ and by
Corollary~\ref{GeneralReductionCorollary} we obtain $f_{q^{''}}(X) -
f_{q^{''}-2}(X) = 2I$ and then by Lemma~\ref{GeneralReduction}
linear growth follows. $\Box$

\bigskip

\subsection{Substituting type $E$}\label{FiniteE}

Similar results can be shown for type $E$. However, all we need for the application to
 representation theory that we have in mind, are some simple facts that an easy calculation by hand can provide.
 Namely, we need to know when  $f_k(E_i )$ is equal to zero for the first time when $i = 6, 7$ or  8 where
  $E_i$ for $i = 6,7,8$ is the adjacency matrix of the Dynkin diagram of type $\E$ as given in section~\ref{Dynkin}.
    This  happens exactly  for $f_{11}(E_6)$, $f_{17}(E_7)$ and $f_{29}(E_8)$.

Furthermore we need the explicit expressions of the matrices
preceding this first zero matrix. They are given by
$$f_{10}(E_6) = ( \varepsilon_6, \varepsilon_5, \varepsilon_3, \varepsilon_4, \varepsilon_2, \varepsilon_1)^T, $$
$$f_{16}(E_7) = I_7, $$
$$ f_{10}(E_8) = I_8 $$
where $ \varepsilon_i$ is as before, that is, it is the row vector that has a one in place $i$ and zero elsewhere.



\bigskip



\subsection{Summary}

 In the preceding sections we have calculated  explicitly
 the Chebyshev polynomials evaluated on the adjacency matrices of the  Dynkin diagrams of types
$A, D, E $, of the extended Dynkin diagrams $\widetilde{A},
\widetilde{D},$ as well as on the diagrams of type $L$, $ \widetilde{L}$ and $\widetilde{DL}$.
 In particular we have shown  that the
Dynkin diagrams and the type $L$ diagram give rise to periodic families and that the
extended Dynkin diagrams and the diagrams of types $\widetilde{L}$ and $\widetilde{DL}$, give rise to families that have linear
growth.

\begin{Theorem} Let $\mathcal D$ be the adjacency matrix of a Dynkin diagram of type $A, D, E$ or a diagrams of type $L$. Then for $r \geq 1$,
we have $f_k({\mathcal D})=0$ if and only if  $k = rh-1$ where $h$ denotes the
Coxeter number of the associated diagram.
\end{Theorem}













\begin{Theorem}\label{q_Extended}
Let $\mathcal D$ be the adjacency matrix of an extended Dynkin
diagram or of a diagram of type $\widetilde{L}$ or $\widetilde{DL}$. Then for $m = rq + u$ where $-1 \leq u \leq q -2$ and $2
\leq r $ we have the recurrence relation

$$ f_{rq +u} (\mathcal D) = r f_{q+u} (\mathcal D) - (r-1) f_{u}(\mathcal
D)$$

(a) if $\mathcal D = \widetilde{A}_n$  then $q=n+1$;

(b) if $\mathcal D = \widetilde{D}_n$ then $q=2n-4$
for $n$ even and $q=8n-16 $ for $n$ odd;

(c) if $\mathcal D = \widetilde{L}_n$ then $q=n$;

(d) if $\mathcal D = \widetilde{DL}_n$ then $q=8n-4$;

(e)     if $\mathcal D = \widetilde{E}_6$ then $q=12$, if $\mathcal
D = \widetilde{E}_7$ then
    $q=72$, and if $\mathcal D = \widetilde{E}_8$  then
     $q=60$.
     \end{Theorem}

{\bf Remark:} It is easily checked that for $\mathcal D$ one of the
extended types $\widetilde{E}$ with $q$ given as above we have $f_{q}
(\mathcal D) - f_{q-2} (\mathcal D) = 2I$.

\section{Application to representation theory}\label{Application}

In this section we show how the previous results can be used to
calculate the minimal projective resolutions of the simple modules
of a class of symmetric algebras -- namely those that are of radical
cube zero and  of tame or of finite representation type. Our method
gives the indecomposable projective components in each degree of the
projective resolutions through a description of the radical layers
of the syzygies.

\bigskip

Let $K$ be a field and let $\Lambda$ be a finite dimensional $K$-algebra.
Then $\Lambda$ is symmetric if there exists a linear map $\nu : \Lambda
\rightarrow k$ such that for all $a,b \in \Lambda$, $\nu(ab) = \nu(ba)$
such that $\rm{Ker}(\nu)$ contains no non-zero left or right ideal.
Recall that the Jacobson radical $J(\Lambda)$ of the algebra
is the smallest ideal of $\Lambda$ such that the quotient is semisimple.

Let $\Lambda$ be a finite dimensional symmetric $K$-algebra such
that $J^3(\Lambda) = 0$. We assume that $\Lambda$ is indecomposable, and
that $J^2(\Lambda) \neq 0$. These algebras are classified in
\cite{Benson} according to the minimal projective resolution of
non-projective finitely generated $\Lambda$-modules. Namely, if
$S_1, S_2, \ldots, S_n$ are the simple $\Lambda$-modules, let
${\mathcal E}_n = (e_{ij})_{i,j= 1, \ldots, n}$ where $e_{ij} = \dim
\;\Ext_\Lambda^1(S_i,S_j)$, which is a symmetric matrix.
Then if the  largest
eigenvalue $\lambda$ of ${\mathcal E}_n$ is $> 2$,
the dimensions in the minimal projective resolutions of the
non-projective finite-dimensional  $\Lambda$-modules are unbounded,
and the algebra is of wild representation type. If  $\lambda = 2$,
the algebra is of tame representation type and the minimal
projective resolutions either are bounded or grow linearly, and if
$\lambda < 2$ they are bounded.
  In the latter two cases the algebras are classified by Dynkin diagrams  $\A, \D, \E$ or the graph
$\TL$ for the finite representation type,   and by the extended Dynkin
diagrams
 $\widetilde{\A}, \widetilde{\D}, \widetilde{\E}$ or the graphs $\widetilde{\TL}$ and
$\widetilde{\bf{\rm DL}}$ for the tame representation type.

\bigskip

To each of these diagrams $\mathcal D$
we associate the quiver $Q({\mathcal D})$, which is
obtained by
replacing
each edge of the
diagram by a pair of  arrows pointing in opposite
directions.
Let $K$ be an algebraically closed field  and let $KQ({\mathcal D})$
be the path algebra of $Q({\mathcal D})$. Then one has an ideal $I$
of   $KQ({\mathcal D})$ such that the corresponding quotient
algebra $\Lambda = KQ({\mathcal D})/I$ is a symmetric algebra such
that $J^3(\Lambda) = 0$.

\bigskip

We recall some properties of indecomposable $\Lambda$-modules as
stated in~\cite{Benson}. The Loewy length
of a $\Lambda$-module is at most  3. If $P$ is a
projective indecomposable $\Lambda$-module, then ${\rm rad}^2(P) = {\rm soc}(P)
\cong  {\rm hd}(P) $.  If $M$ is indecomposable but  not simple or projective,
then the Loewy length of $M$ is equal to 2 and we
denote by $d(M)$ the column vector $\left(
\begin{array}{c}  \alpha (M) \\ \beta (M) \end{array}\right)$ where
$\alpha (M) = \underline{\dim} \; \hd(M)$ and $\beta (M) =
\underline{\dim}  \; \rad(M)$ where $\underline{\dim} V$ denotes the
dimension vector of the finite dimensional $\Lambda$-module $V$.
This is the column vector of length $n$ whose $i$th entry
corresponds to the multiplicity of the simple module $S_i$
as a composition factor in $V$. If $S$ is a simple
$\Lambda$-module then $\alpha (S) = \underline{\dim} \; S$ and $\beta
(S) = 0$. And if $P(S)$ is a projective cover of $S$ then the
dimension vector of its three radical layers are described by the
vector $ \left(
\begin{array}{c} \alpha (S) \\  {\mathcal E}_n\alpha (S) \\ \alpha
(S)
\end{array}\right)$. More generally, the projective cover $P(M)$ of
a non-projective indecomposable $\Lambda$-module $M$ has  in its
radical layers the dimension vectors $\alpha (M)$, ${\mathcal E}_n
\alpha (M)$ and $\alpha (M)$. Therefore if $\Omega(M)$ is non-simple
it has in its radical layers the dimension vectors ${\mathcal E}_n
\alpha (M) - \beta (M)$ and $\alpha (M)$. Furthermore, $d(\Omega
(M)) = B
d(M)$ with  $B = \left( \begin{array}{cc} {\mathcal E}_n & -I_n \\
I_n & 0
\end{array} \right)$ where $I_n$ is the identity $n \times n$
matrix. Therefore for $M$ an indecomposable non-projective
$\Lambda$-module we get the recurrence relation
$$ d(\Omega^k (M)) = B^k d(M) $$
if none of the $\Omega^j(M)$ for $j=0, \ldots, k$ are simple
$\Lambda$-modules. Furthermore, we observe that the entries of the
matrix $B^k$  can be given in terms of Chebyshev polynomials.

\begin{Lemma}~\label{Chebyshev}
Let $M$ be an $n \times n$ matrix and let $B = \left(
\begin{array}{cc} M & -I_n \\
I_n & 0 \end{array} \right)$. Then for $i \geq 0$,
$$B^i = \left( \begin{array}{cc} f_i(M) & -f_{i-1}(M) \\
f_{i-1}(M) & -f_{i-2}(M)   \end{array} \right)$$ where $f_i(x)$ is
the  Chebyshev polynomial of the  second-kind defined by the
recurrence relation $f_i(x) = x f_{i-1}(x) - f_{i-2}(x)$ and with
initial conditions $f_0(x)=1$ and $f_1(x)=x$.
\end{Lemma}

{\it Proof:} By definition of $B$ we have  $f_0(M)=I_n$ and
$f_1(M)=M$ and we pose $f_{-1}(M)=0$. Suppose now the result holds
for $i$. Then
$$B^iB= \left( \begin{array}{cc} f_i(M) & -f_{i-1}(M) \\ f_{i-1}(M) & -f_{i-2}(M)   \end{array} \right)
 \left( \begin{array}{cc} M & I_n \\ -I_n & 0   \end{array} \right)= \left( \begin{array}{cc} f_{i+1}(M) & -f_{i}(M) \\  f_{i}(M) & -f_{i-1}(M)   \end{array} \right). \;\;\;\; \Box$$

\bigskip

If $\Omega (M)$ is simple then $ \Omega (M) = soc (P(M))$ and
$d(\Omega (M)) = (\alpha (M), 0)^T$, however $B d(M) = (0, \alpha
(M))^T$. Thus we know that $\Omega^k(M)$ is simple if the first $n$
entries of the vector $B ^k d(M)$ are zero and if there is only one
non-zero entry in the $n+1$ to $2n$ components  of $B ^k d(M)$.

\bigskip

Let $\Lambda$ be one of the algebras defined above. We say that
$\Lambda$ is of type $\A, \D, \E, \TL, \widetilde{\A},
\widetilde{\D}  $ etc. (or $\A_n, \D_n, \E_6, \E_7, \E_8, \TL_n,
\widetilde{\A}_n, \widetilde{\D}_n $ etc.  if we need to specify the
number of vertices of the diagram) if the underlying  diagram
of the quiver of $\Lambda$ is of that type.

\bigskip

The following two theorems are an application of the results of
section~\ref{Evaluation} . In the case of  Theorem~\ref{period} this gives a proof
of a result that can be deduced from Theorem 2.1 in~\cite{BBK}.

\begin{Theorem}\label{period}
The minimal projective resolutions of the simple $\Lambda$-modules
are periodic

(a) of period $2n$ if $\Lambda$ is of type $\A_n$,

(b) of period $2n$ if $\Lambda$ is of type $\TL_n$,

(c) of period $2n-3$ if $\Lambda$ is of type $\D_n$ for $n$ even and
of period $2(2n-3)$ for $n$ odd,

(d) of period $22$ if $\Lambda$ is of type $\E_6$, of period $17$ if
$\Lambda$ is of type $\E_7$ and of period $29$ if $\Lambda$ is of
type $\E_8$.

\end{Theorem}

\bigskip


\begin{Definition} We call a projective resolution ${\mathbf R^\bullet}$ of a
$\Lambda$-module $M$ of {\em linear growth of factor p} if given the
first $p$ terms
$$R_{p-1} \rightarrow \cdots \rightarrow R_1 \rightarrow R_0 \rightarrow M \rightarrow 0
$$
of ${\mathbf R^\bullet}$, all other terms  of ${\mathbf R^\bullet}$ have the form
$$ R_{kp+l} = N_l^{\oplus k} \oplus R_l \;\; \mbox{for} \;\; 0 \leq l \leq p-1
$$
where $N_l$ is a direct sum of components $R_i$ for $0 \leq i \leq
p$.
\end{Definition}

Let $q$ be as given in Theorem~\ref{q_Extended}.

\begin{Theorem}\label{periodextended}
Let $\Lambda$ be of extended Dynkin type or of type
$\widetilde{DL}$. Then the minimal projective resolutions of the
simple $\Lambda$-modules are of linear growth of factor $q$.
\end{Theorem}

\bigskip

Theorem~\ref{period} follows directly from the calculation of the
syzygies in the next Proposition.

\begin{Proposition}~\label{OmegatypeA}
(a) Let $\Lambda$ be of type $\A_n$. Then for all simple
$\Lambda$-modules $S_i$ we have $\Omega^j (S_i)$ is not simple for $j < n$ and $\Omega^n (S_i) = S_{n-i+1}$.

(b) Let $\Lambda$ be of type $\TL_n$. Then for all simple
$\Lambda$-modules $S_i$, we have $\Omega^j (S_i)$ is not simple for $j < 2n$ and $\Omega^{2n} (S_i) = S_{i}$.

(c) Let $\Lambda$ be of type $\D_n$. Then for all simple
$\Lambda$-modules $S_i$, we have  $\Omega^j (S_i)$ is not simple for $j < 2n-3$ and
$$ \Omega^{2n-3} (S_i) = \left\{ \begin{array} {ll}
S_2 & \mbox{if $i=1$ and $n$ is odd} \\ \\
S_1 & \mbox{if $i=2$ and $n$ is odd} \\ \\
S_i & \mbox{otherwise.} \\
\end{array} \right.
$$

(d) Let $\Lambda$ be of type $\E_6$ , $\E_7$ or $\E_8$. Then for all
simple $\Lambda$-modules $S_i$, we have
 $$ \Omega^{j} (S_i) = \left\{ \begin{array} {ll}
S_i & \mbox{if $n =6$ and $j= 11$ and for $i = 3,4$} \\ \\
S_{n-i+1} & \mbox{if $n =6$ and $j= 11$ and for $i =  1,2,5,6$} \\ \\
S_i & \mbox{if $n =7$ and $j= 17$ and for all $i$} \\ \\
S_i & \mbox{if $n =8$ and $j= 29$ and for all $i.$} \\ \\
\end{array} \right.
$$
\end{Proposition}

{\it Proof:} The Proposition follows directly from
Lemma~\ref{Chebyshev} and the results of section~\ref{Evaluation}
$\Box$.

\bigskip







{\it Proof of Theorem~\ref{periodextended}:} It is enough to
establish the recurrence formula for the Chebyshev polynomials
evaluated on the adjacency matrices in question, since by
Lemma~\ref{Chebyshev} it is then straightforward to determine the
projective resolutions. Namely, if the algebra has $n$ simple
modules, then the indecomposable projectives constituting the $r$th
term in a projective resolution of the simple $S_i$ are exactly
given by the first $n$ entries of the  $i$th row of the matrix
$B^r$.

 Let ${\mathcal D}$ be the adjacency matrix of the graph
underlying $\Lambda$ (recall that ${\mathcal D}$ is not of type
$\widetilde{L}$). Then it follows from Theorem~\ref{q_Extended} in
combination with Lemma~\ref{GeneralReduction} that given
$f_0({\mathcal D}), \ldots, f_q({\mathcal D})$ we have for $q+1 \leq
n \leq 2q-1$, where $n=q+l$ for $1 \leq l \leq q-1$, that
$f_{q+l}({\mathcal D}) = (f_{q-2-l}({\mathcal D}) + f_l({\mathcal
D})) + f_l({\mathcal D})$ since in this case $Z=2I$. So if we
calculate a minimal projective resolution of the simple $S_i$, then
the components of $q+l$th term, for $1 \leq l \leq q-1$,  are
$R_{q+l} = N_l \oplus R_l$ where $R_l$ is determined by the $i$th
row of $f_l({\mathcal D})$ and $N_l = R_{q-2-l} \oplus R_l$ is
determined by the $i$th row of $f_{q-2-l}({\mathcal D}) +
f_l({\mathcal D})$. For $n \geq 2q-1$ we have the following: Write
$n = rq+l$ with $-1 \leq l \leq q-2$. Furthermore, note that here
$Z=2I = 2f_0({\mathcal D})$. Then by Lemma~\ref{GeneralReduction}
and the remark that follows we have, for $-1 \leq l \leq q-2$, that
$R_{rq+l} = N^{\oplus r}_l \oplus R_l$ where $N_l = R_{q-2-l} \oplus
R_l$.
$\Box$

\bigskip

{\bf Remark:} Suppose $\Lambda$ is of type
$\widetilde{L}$. Then the projective resolution of the simple
$\Lambda$-modules is almost of linear growth. More precisely in that
case we have by Corollary~\ref{extendedL} that $f_{q}(X)- f_{q-2}(X)
= Z$ where $Z=2I^0 = 2f_0({\mathcal D})^0$. Given the first $q$
terms in a projective resolution of a simple $\Lambda$-module $S$ we
then obtain by Lemma~\ref{GeneralReductionZero} and the remark that
follows the $n$th term in the projective resolution of $S$, for
$n\geq 2q-1$, in the following way: if we write $n = rq+l$ with $-1
\leq l \leq q-2$, then
$$
 N_{rq+l} = \left\{\begin{array}{ll}
(N_{l}^0)^{\oplus r} + R_l
 & r \mbox{ even } \cr
N_{l}^{\oplus r} + R_l^0 & r \mbox{ odd }\end{array}\right.
$$
where $N_l$ is as defined in the proof above and the components in
$N_l^0$ are given by the $i$th row of $I^0(f_{q-2-l}({\mathcal D})+
f_l({\mathcal D}))I^0$ and $R_l^0$ is given by the $i$th row of
$I^0f_l({\mathcal D})I^0$.


\bigskip



\section{Evaluating Chebyshev polynomials on positive symmetric
matrices}\label{Symmetric}

We conclude by a general statement on evaluating the Chebyshev
polynomials $(f_k(x))_k$ on positive symmetric matrices.

\begin{Theorem}
Assume $X$ is a symmetric matrix, with entries in $\mathbb{Z}_{\geq 0}$
and assume $X$ is indecomposable. Then

a) $f_d(X)=0$ for some $d\geq 1$ if and only if $X$ is the adjacency
matrix of a diagram of type A, D, E or L.

b) The family $(f_k(X))_k$ grows linearly if and only if
$X$ is the adjacency matrix of a diagram of type $\widetilde{A}, \widetilde{D}, \widetilde{E}, \widetilde{L}$ or $\widetilde{DL}$.
\end{Theorem}

{\it Proof:} a) If $X$ is the adjacency matrix of a Dynkin diagram
or of a graph of type $\TL$, then the result follows from
section~\ref{Evaluation}. Suppose that $f_d(X) = 0 $ for some $n$,
then by Lemma~\ref{Periodic} $X$ annihilates the sequence of
polynomials $f_k(x)$ periodically. Since $X$ is a symmetric integer
matrix with entries in $\mathbb{Z}_{\geq 0}$, then it is the
adjacency matrix of a finite connected graph. To this graph we can
associate a unique symmetric algebra $\Lambda$ with radical cube
zero, such that the sequence of polynomials $(f_k(X))_{k\geq 0}$
describes the growth of a minimal projective resolution of the
simple modules, as explained at the beginning of section 4. To
construct this algebra, we replace each edge by a pair of arrows
pointing in opposite directions. Then the algebra is the path
algebra modulo the ideal generated by quadratic relations in the
arrows, and there is a unique choice of such relations making the
algebra symmetric with radical cube zero.

Then it follows
 from
section~\ref{Application} that for such an algebra the minimal
projective resolutions of the simple modules are periodic  and thus
following \cite[1.1]{Benson} $\Lambda$ is of type $\A, \D, \E$ or
$\TL$.

b) If $X$ is of extended Dynkin type or of types $\widetilde{L}$ or
$\widetilde{DL}$ then by section~\ref{Evaluation}, $(f_k(X))_k$
grows linearly. Conversely, suppose that $(f_k(X))_k$ grows
linearly. As in a), $X$ gives rise to a symmetric algebra $\Lambda$
of radical cube zero. By section~\ref{Application}  this implies
that the minimal projective resolutions of the simple
$\Lambda$-modules grow linearly. Following \cite[1.1]{Benson}
$\Lambda$ is of type $\widetilde{\A}, \widetilde{\D},
\widetilde{\E}$, $\widetilde{\TL}$, or $\widetilde{\bf{\rm DL}}$.
$\Box$

\section*{Appendix}

 Dynkin diagrams with labels:

{\tiny

$$\begin{array}{ll} A_n, n\geq 2 :   \xymatrix{
\overset{1}{\bullet} \ar@{-}[r]  & \overset{2}{\bullet} \;\;\;\;\;
\cdots \;\;\;\;\; \stackrel{n-1}{\bullet} \ar@{-}[r] &
\stackrel{n}{\bullet} } \;\;\; & \;\;\;
 L_n, n \geq 2 :   \;\;\;\;\;\;\;\;
\xymatrix{ \overset{1}{\bullet} \ar@(ul,dl)@{-} \ar@{-}[r]  &
\overset{2}{\bullet} \;\;\;\;\; \cdots \;\;\;\;\;
\stackrel{n-1}{\bullet} \ar@{-}[r] & \stackrel{n}{\bullet} }
\\ & \\ & \\ & \\
D_n, n\geq 4 :   \xymatrix{
\overset{1}{\bullet} \ar@{-}[dr] \\
 &  \overset{3}{\bullet}  \ar@{-}[r] & \overset{4}{\bullet} \;\;\;\;\;
\cdots \;\;\;\;\; \stackrel{n-1}{\bullet} \ar@{-}[r] &
\stackrel{n}{\bullet} \\
\overset{2}{\bullet} \ar@{-}[ur] } \;\;\; & \;\;\; E_6 : \xymatrix{
&& \overset{6}{\bullet} \ar@{-}[d] \\
 \overset{1}{\bullet} \ar@{-}[r]  & \overset{2}{\bullet}
\ar@{-}[r] &   \overset{3}{\bullet}  \ar@{-}[r] &
\stackrel{4}{\bullet} \ar@{-}[r] & \stackrel{5}{\bullet} }
\\
E_7 :   \xymatrix{
&& \overset{7}{\bullet} \ar@{-}[d] \\
 \overset{1}{\bullet} \ar@{-}[r]  & \overset{2}{\bullet}
\ar@{-}[r] &   \overset{3}{\bullet}  \ar@{-}[r] &
\stackrel{4}{\bullet} \ar@{-}[r] & \stackrel{5}{\bullet}  \ar@{-}[r]
& \stackrel{6}{\bullet}} \;\;\; & \;\;\; E_8 :   \xymatrix{
&& \overset{8}{\bullet} \ar@{-}[d] \\
 \overset{1}{\bullet} \ar@{-}[r]  & \overset{2}{\bullet}
\ar@{-}[r] &   \overset{3}{\bullet}  \ar@{-}[r] &
\stackrel{4}{\bullet} \ar@{-}[r] & \stackrel{5}{\bullet} \ar@{-}[r]
& \stackrel{6}{\bullet} \ar@{-}[r] & \stackrel{7}{\bullet}}
\end{array}$$

}

Extended Dynkin diagrams with labels:

{\tiny
$$ \begin{array}{ll}\widetilde{A}_n, n\geq 2 :   \xymatrix{
&   \overset{n+1}{\bullet} \ar@{-}[dl] \ar@{-}[dr]\\
 \overset{1}{\bullet} \ar@{-}[r]  & \overset{2}{\bullet}
\;\;\;\;\; \cdots \;\;\;\;\; \stackrel{n-1}{\bullet} \ar@{-}[r] &
\stackrel{n}{\bullet} } \;\;\; & \;\;\; \widetilde{L}_n, n \geq 2 :
\;\;\;\;\;\;\;\;  \xymatrix{ \overset{1}{\bullet} \ar@(ul,dl)@{-}[]
\ar@{-}[r]  & \overset{2}{\bullet} \;\;\;\;\; \cdots \;\;\;\;\;
\stackrel{n-1}{\bullet} \ar@{-}[r] & \stackrel{n}{\bullet}
\ar@(ur,dr)@{-}[]  } \\ & \\ & \\
\widetilde{D}_n, n\geq 3 : \xymatrix{
\overset{1}{\bullet} \ar@{-}[dr] &&&& \overset{n}{\bullet} \ar@{-}[dl]\\
 &  \overset{3}{\bullet}  \ar@{-}[r] & \overset{4}{\bullet} \;\;\;\;\;
\cdots \;\;\;\;\; \stackrel{n-2}{\bullet} \ar@{-}[r] &
\stackrel{n-1}{\bullet} \\
\overset{2}{\bullet} \ar@{-}[ur] &&&& \overset{n+1}{\bullet}
\ar@{-}[ul] } \;\;\; & \;\;\; \widetilde{DL}_n, n\geq 3 : \xymatrix{
\overset{1}{\bullet} \ar@{-}[dr] \\
 &  \overset{3}{\bullet}  \ar@{-}[r] & \overset{4}{\bullet} \;\;\;\;\;
\cdots \;\;\;\;\; \stackrel{n-1}{\bullet} \ar@{-}[r] &
\stackrel{n}{\bullet} \ar@(ur,dr)@{-}[]\\
\overset{2}{\bullet} \ar@{-}[ur]   }\\ & \\ & \\
\widetilde{E}_6 : \xymatrix{
&& \overset{7}{\bullet} \ar@{-}[d] \\
&& \overset{6}{\bullet} \ar@{-}[d] \\
 \overset{1}{\bullet} \ar@{-}[r]  & \overset{2}{\bullet}
\ar@{-}[r] &   \overset{3}{\bullet}  \ar@{-}[r] &
\stackrel{4}{\bullet} \ar@{-}[r] & \stackrel{5}{\bullet} }\;\;\; &
\;\;\; \widetilde{E}_7 :   \xymatrix{
&&& \overset{7}{\bullet} \ar@{-}[d] \\
 \overset{8}{\bullet} \ar@{-}[r] &\overset{1}{\bullet} \ar@{-}[r]  & \overset{2}{\bullet}
\ar@{-}[r] &   \overset{3}{\bullet}  \ar@{-}[r] &
\stackrel{4}{\bullet} \ar@{-}[r] & \stackrel{5}{\bullet}  \ar@{-}[r]
& \stackrel{6}{\bullet}}\\ & \\ & \\ & \\
\widetilde{E}_8 :   \xymatrix{
&& \overset{8}{\bullet} \ar@{-}[d] \\
 \overset{1}{\bullet} \ar@{-}[r]  & \overset{2}{\bullet}
\ar@{-}[r] &   \overset{3}{\bullet}  \ar@{-}[r] &
\stackrel{4}{\bullet} \ar@{-}[r] & \stackrel{5}{\bullet} \ar@{-}[r]
& \stackrel{6}{\bullet} \ar@{-}[r] & \stackrel{7}{\bullet}
\ar@{-}[r] & \stackrel{9}{\bullet}} \end{array}$$ }

\end{document}